\documentclass[11pt,
]{article}

\usepackage{amssymb,amsfonts,amstext,amsmath,amsthm,
latexsym,mathrsfs,amsbsy}
\usepackage{bbold}

\usepackage{marginnote}

\usepackage{mparhack}
\usepackage{esvect}

\input{47.sty}

\vyk{
\usepackage{graphicx}
\usepackage{xcolor}

\usepackage
[%
hypertex
,citecolor=blue
]
{hyperref}
}

\begin{document}

\title
{On intermediate extensions of generic extensions 
by a random real}

\author 
{
Vladimir Kanovei\thanks{IITP RAS and MIIT,
  Moscow, Russia, \ {\tt kanovei@googlemail.com} --- 
contact author. 
Partial support of   RFFI grant 17-01-00705 acknowledged.
}  
\and
Vassily Lyubetsky\thanks{IITP RAS,
  Moscow, Russia, \ {\tt lyubetsk@iitp.ru}. 
Partial support of Russian Scientific Fund grant 14-50-00150 
acknowledged. 
}
}

\date 
{\today}

\maketitle

\begin{abstract}
The paper is the second of our series of notes aimed 
to bring back in circulation some bright ideas of 
early modern set theory, mainly due to Harrington 
and Sami, which have never been adequately presented 
in set theoretic publications.
We prove that if a real $a$ is random over a model $M$ 
and $x\in M[a]$ is another real then either (1) $x\in M$, 
or (2) $M[x]=M[a]$, or (3)
$M[x]$ is a random extension of $M$ and $M[a]$ 
is a random extension of $M[x]$. 
This is a less-known result of old 
set theoretic folklore, and, as far as we know, 
has never been published. 

As a corollary, we prove that $\fs1n$-Reduction 
holds for all $n\ge3$, in a model extending $\rL$ by 
$\ali$-many random reals.


{\footnotesize
\tableofcontents 
}

\end{abstract}

\np

\parf{Introduction}

It is known from Solovay \cite{sol}, and especially 
Grigorieff \cite{gri} in most general form, that any 
subextension $\rV[x]$ of a generic extension 
$\rV[G]$, generated by a set $x\in\rV[G]$, 
is itself a generic extension $\rV[x]=\rV[G_0]$ 
of the same ground universe $\rV$, 
and the whole extension $\rV[G]$ is equal to a 
generic extension $\rV[G_0][G_1]$ of the 
subextension $\rV[x]=\rV[G_0]$. 
See a more recent treatment of this question in 
\cite{jechmill,zapt,kl21,kl32e}. 
In particular, it is demonstrated in \cite{kl32e} 
that if $\dP=\stk\dP\leq\in \rV$ is a forcing notion,  
a set $G\sq\dP$ is $\dP$-generic over $\rV$, 
$t\in\rV[G]$ is a $\dP$-name, 
$x=t[G]\in\rV[G]$ is the $G$-valuation of $t$, 
and $x\sq\rV$, then 
\ben
\Renu
\itlb{i1}%
there is a set $\Sg\sq\dP$ such that 
$\rV[\Sg]=\rV[x]$  and 
$G$ is $\Sg$-generic over $\rV[x]$;

\itlb{i2}%
there exists a stronger order $\leq_t$ on $\dP$ 
(so that $p\leq q$ implies $p\leq_t q$) 
in $\rV$ such that $\Sg$ itself is 
$\stk\dP{\leq_t}$-generic over $\rV[\Sg]=\rV[x]$. 
\een  
However the nature and forcing properties  
of the derived forcing notions 
$\dP_0=\stk\dP{\leq_t}\in\rV$ 
and $\dP_1(x)=\stk\Sg{\leq}\in\rV[x]$ 
is not immediately clear. 

At the trivial side, we have the Cohen forcing 
$\dP=\dC=\bse.$
In this case, $\dP_0$ and $\dP_1(x)$ are 
countable forcing notions, hence the corresponding 
extensions, $\rV\to\rV[x]$ and $\rV[x]\to \rV[G]$ 
in the above scheme, 
are Cohen generic or trivial. 
As observed in \cite{kl32e}, this leads to the 
following result of set theoretic folklore, never 
explicitly appeared in set theoretic publications, 
except for 
\cite[Lemma 1.9]{samiPHD}. 
(It can also be derived from some results in \cite{gri}, 
especially 4.7.1 and 2.14.1.) 

\bte
[folklore, Sami]
\lam{intC}
Let\/ $a\in\dn$ be Cohen-generic over the ground 
set universe\/ $\rV$. 
Let\/ $x$ be a real in\/ $\rV[a]$. 
Then we have exactly one of the following$:$
\ben
\tenu{{\rm(C\arabic{enumi})}}
\itlb{intC1}%
$x\in\rV\,;$ \qquad\qquad
{\rm(C2)} $\,\rV[x]=\rV[a]\,;$

\atc
\itlb{intC2}%
{\rm (a)} 
$\rV[x]$ is a Cohen-generic extension of\/ $\rV$,   
and\/ \\[0.5ex] 
{\rm (b)} 
$\rV[a]$ 
is a Cohen-generic extension of\/ $\rV[x]$.\snos
{Theorem~\ref{intC} 
dramatically fails for intermediate extensions 
not generated by sets, \cite{abyss}.}
\qed 
\een
\ete

A much more complex case is the Levy -- Solovay 
extension of $\rL$, the constructible universe.  
As established in \cite{sol}, such an extension is 
equal to a Levy -- Solovay extension of $\rL[x]$ 
for any real $x$ it contains.

The following theorem, proved below, 
is a result of the same type. 

\bte
\lam{intR}
Let\/ $a\in\dn$ be Solovay-random over the ground 
set universe\/ $\rV$, 
Let\/ $x$ be a real in\/ $\rV[a]$. 
Then we have exactly one of the following$:$
\ben
\tenu{{\rm(R\arabic{enumi})}}
\itlb{intR1}%
$x\in\rV\,;$ \qquad\qquad
{\rm(R2)} $\,\rV[x]=\rV[a]\,;$

\atc
\itlb{intR2}%
{\rm (a)} 
$\rV[x]$ is a Solovay-random extension of\/ $\rV$,   
and\/ \\[0.5ex] 
{\rm (b)}    
$\rV[a]$ is a Solovay-random extension of\/ 
$\rV[x]$. 
\een
\ete

It is {\bfit not\/} asserted though that the 
real $x$ itself is 
random over $\rV$ in (a) and/or the real 
$a$ itself is random over $\rV[x]$ in (b).
  
Note that Theorem~\ref{intR} contains two separate 
dichotomies: \ref{intR1} vs.\ \ref{intR2}(a) and 
(R2) vs.\ \ref{intR2}(b).
In spite of obvious semblance of Theorem~\ref{intC}, 
this theorem takes more effort. 
Its proof (it begins in Section~\ref{XX1}) involves 
some results related rather to real analysis 
and measure theory. 

Now we proceed with an application of Theorem~\ref{intR}.

\parf{A corollary: 
Reduction in extensions by random reals}
\las{redf}

\vyk{
The separation property for a pointclass $K$, or 
simply \rit{\dd KSeparation}, is 
the assertion that any two disjoint sets $X,Y$ 
in $K$ (in the same Polish space) can be separated 
by a set in $K\cap \dop K$, where $\dop K$ is the 
pointclass of complements of sets in $K$. 
}%

The reduction property for a pointclass $K$, or 
simply \rit{\dd KReduction}, is 
the assertion that for any two sets $X,Y$ in $K$ 
(in the same Polish space) there exist 
\rit{disjoint} sets $X'\sq X$, $Y'\sq Y$  
in the same class $K$, such that $X'\cup Y'=X\cup Y$. 

It is known classically from 
studies of 
Kuratowski \cite{kursep} that 
Reduction holds for  
$\fp11$ and $\fs12$, but fails for 
$\fs11$ and $\fp12$. 
%
As for the higher projective classes, 
Addison \cite{add2} proved that the axiom 
of constructibility $\rV=\rL$ implies that 
Reduction holds for  
$\fs1n\yd n\ge 3$, but fails for 
$\fp1n\yd n\ge 3$. 
On the other hand, by Martin \cite{martAD}, 
the axiom of projective determinacy $\mathbf{PD}$ 
implies that, similarly to projective level $1$, 
$\fp1n$-Reduction holds for all odd numbers $n\ge3$, 
and, similarly to projective level $2$, 
$\fs1n$-Reduction holds for all even $n\ge4$. 

Apparently not much is known on 
Reduction for higher projective classes in generic 
models. 
One can expect that rather homogeneous, well-behaved 
forcing notions produce generic extensions of $\rL$, 
in which Reduction keeps to be true for projective 
classes $\fs1n$ and accordingly fails for $\fp1n$, 
$n\ge3$, while in specially designed non-homogeneous 
extensions this pattern can be violated. 
This idea is supported by a few known results. 
Ramez Sami \cite{samiPHD} proved 

\bte
[Sami]
\lam{mtC}
It is true in any extension of\/ $\rL$ by\/ $\ali$  
Cohen reals that if\/ $n\ge3$ then\/ 
$\is1n$-Reduction holds, and hence\/ 
$\fs1n$-Reduction holds, too.\snos
{To prove that $\is1n$-Reduction implies  
the boldface $\fs1n$-Reduction, it suffices to use  
a double-universal pair of $\is1n$ sets, as those 
used in a typical proof that $\fs1n$-Reduction 
and $\fs1n$-Separation contadict each other.  
This argument does not  work 
for Separation though.}\qed
\ete

On the other hand, we proved in \cite{kl28}
that Reduction fails for $\fs13$ 
(and in fact Separation fails for both $\fs13$ 
and $\fp13$)
in a rather complicated model related to an
$\ali$-product of forcings similar to 
Jensen's minimal forcing \cite{jenmin}.
See also \cite{kl30e,kl38} on similar models in which 
the Uniformization principle fails for $\fp12$ 
(or $\fp1n$ for a given $n\ge3$) sets 
with countable sections, and \cite{kl47a} on some 
related (and very complex) models of Harrington.
Here we prove the following theorem.      

\bte
\lam{mt'}
It is true in any extension of\/ $\rL$ by\/ $\ali$  
Solovay-random reals that if\/ $n\ge3$ then\/ 
$\is1n$-Reduction holds, and hence\/ 
$\fs1n$-Reduction holds, too.
\ete

Note that the theorem also holds in models obtained 
by adding any uncountable (not necessarily $\ali$) 
number $\ka$ of random reals. 
(Because such models are elementarily equivalent to the 
extension by $\ali$ random reals.)

Sami's proof of Theorem~\ref{mtC} involves 
Theorem~\ref{intC}.
Accordingly, we'll use Theorem~\ref{intR} 
rather similar way.
The following lemma is the key ingredient.

\ble
[proof see Section~\ref{lok}]
\lam{llok}
If\/ $n\ge 2$ and\/ $\vpi(x)$ is a parameter-free\/ 
$\is1n$ formula then there is a parameter-free\/ 
$\is1n$ formula\/ $\vpa(x)$ such that if\/ $x$ is 
a real in an\/ $\ali$-random extension\/ $N$ 
of\/ $\rL$ then\/ $\vpi(x)$ holds in\/ $N$ iff\/ 
$\rL[x]\mo\vpa(x)$. 
\ele

A similar result was obtained by Solovay \cite{sol} 
(\poo\ Levy -- Solovay extensions) and by 
Sami \cite{samiPHD} 
(\poo\ extensions by $\ali$ Cohen reals). 

\bpf[Theorem \ref{mt'}, sketch]
The idea, due to Sami \cite[Lemma 1.11]{samiPHD}, 
is to closely emulate 
Addison's proof of $\is1n$-Reduction in $\rL$. 

Arguing in an\/ $\ali$-random extension\/ $N$ 
of\/ $\rL$, we suppose that $n\ge3$, and 
$X=\ens{x}{\vpi(x)}$ and $Y=\ens{x}{\psi(x)}$
are sets of reals, $\vpi$ and $\psi$ being $\is1n$ 
formulas. 
Then, by Lemma~\ref{llok}, we have 
$X=\ens{x}{\rL[x]\mo\vpa(x)}$ and 
$Y=\ens{x}{\rL[x]\mo\psa(x)}$, where 
$\vpa$ and $\psa$ are still $\is1n$-formulas. 
Thus $\vpa(x)$ is $\sus y\,\Phi(x,y)$ and 
$\psa(x)$ is $\sus y\,\Psi(x,y)$, 
$\Phi$ and $\Psi$ being $\ip1{n-1}$.

Still arguing in $N$, 
if $x\in\dn$ then let $\llx x$ be the 
canonical G\"odel wellordering of the reals in 
$\rL[x]$, of order type $\omi$. 
The crucial property of this system of order relations 
says that the \rit{bounded quantifiers}
$\kaz y'\llx x y$ and $\kaz y'\lelx x y$, applied to 
a $\is1n$ formula, yield a $\is1n$ formula. 
It follows that the sets
$$
\bay{rcl}
X' &=& \ens{x}{\rL[x]\mo \sus y
\big(
\Phi(x,y)\land \kaz y'\llx x y \:\neg\:\Psi(x,y')
\big)}\\[1ex]
Y' &=& \ens{x}{\rL[x]\mo \sus y
\big(
\Psi(x,y)\land \kaz y'\lelx x y \:\neg\:\Phi(x,y')
\big)}
\eay
$$
are $\is1n$, because the relativization to $\rL[x]$ 
does not violate being $\is1n$ ($n\ge 2$).\vtm%

\epF{Theorem~\ref{mt'}, modulo Lemma~\ref{llok} 
and Theorem~\ref{intR}}

\parf{Randomness is measure-independent}
\las{XX0}  

\rit{Random} 
(or \rit{Solovay-random}) reals, 
over a set universe $\rV$, 
are usually defined as those reals in $\dn,$ 
or true reals in the unit interval $[0,1]=\dI$, 
which avoid Borel sets, coded in $\rV$ and null 
with respect to, resp., 
the usual product probability measure $\mu$ on $\dn,$ 
or the Lebesgue measure $\la$ on $\dI$. 

That the $\mu$-random reals in $\dn$ and $\la$-random 
reals in $\dI$ produce the same generic extensions 
and thereby both notions can be identified, is 
witnessed by the Borel map 
$f(a)=\sum_{a(n)=1}2^{-n-1}:\dn\onto\dI$. 
It satisfies $\la(\im fX)=\mu(X)$ for any Borel 
$X\sq\dn,$ therefore if $a\in\dn$ and $x=f(a)\in\dI$ 
then $a$ is $\mu$-random iff $x$ is $\la$-random, 
and $\rV[a]=\rV[x]$, of course. 
There is a general version of such a correspondence, 
which will be used in the proof of 
Theorem~\ref{intR} below.

\ble
\lam{xxl}
Assume that\/ $\nu$ is a continuous\/ 
{\rm(that is, all singletons are null sets)} 
Borel probability measure 
defined on\/ $\dn$ in a set universe\/ $\rV$. 
Then there is a Borel map\/ $g:\dn\onto\dI$, 
coded in\/ $\rV$, and such that if\/ $a\in\dn$ 
and\/ $x=g(a)\in\dI$ then\/ $a$ is $\nu$-random 
over\/ $\rV$ iff\/ $x$ is\/ $\la$-random 
over\/ $\rV$, and\/ $\rV[a]=\rV[x]$. 
\ele
\bpf
Let $\lx$ be the lexicographical order on $\dn,$ 
and let $\inl ab=\ens{a'}{a\lx a' \lx b}$ denote 
$\lx$-intervals. 
Let $g(a)=\mu(\inl{\olex}{a})$, 
where $\olex\in\dn$ is the $\lx$-least element, 
$\olex(k)=0\yt\kaz k$.
Easily $g$ is measure-presirving: if $X\sq\dn$ 
is Borel then $\nu(X)=\la(\im gX)$. 
(See \eg\ the proof of Theorem 17.41 in 
Kechris~\cite{dst}.) 
It follows that $a$ is $\nu$-random 
iff $x$ is $\la$-random, 
whenever $a\in\dn$ and $x=g(a)$. 
To see that $a\in\rV[x]$, note that $J=\aim gx$ 
is a closed $\lex$-interval in $\dn$, the interior 
of which (if non-empty) is a $\nu$-null set, 
hence $a$ is equal to oneof the two endpoints 
of $J$.
\epf

\parf{Intermediate submodels of random extensions:
kase split}
\las{XX1}  

We begin here a proof 
of Theorem~\ref{intR}. 
It will use only basic forcing ideas and some 
classical theorems related to real analysis.  

Thus let $\jao\in\dn$ be Solovay-random over the 
background set universe $\rV$. 
We shall assume that ${\jxo}\in\rV[\jao]$ is a real in the 
unit segment $[0,1]$ of the true real line $\dR$. 
As the Solovay-random forcing admits continuous reading 
of names, there is a continuous map 
$f:\dn\to \dI$, coded in $\rV$, 
such that ${\jxo}=f({\jao})$. 
Let $\muo$ be the usual product probability measure 
on $\dn,$ 
and $\la$ be the Lebesgue measure on the 
segment $\dI=[0,1]$. 

{\ubf We have to prove the trichotomy 
$\ref{intR1}$ vs.\ {\rm(R2)} vs.\ $\ref{intR2}$ 
of Theorem~\ref{intR}}.\vom

{\ubf First split.} 
\rit{Arguing in\/ $\rV$}, consider the 
set $C=\ens{x\in\dI}{\muo(\aim fx)>0}$. 
It is at most countable.  
Consider the complementary sets 
$D=\aim fC$ and $A_1=\dn\bez D$.
These are resp.\ $\Fs$ and $\Gd$ sets coded in $\rV$, 
we identify them with \lap{the same} 
(\ie, coded by the same codes) 
sets in the extensions $\rV[\jao]$, $\rV[\jxo]$.

{\ubf Case 1:} 
$\jao\in D$. 
Then there is a real $\bar y\in\dI\cap\rV$ 
such that $\jao\in\aim f{\bar y}$, 
hence ${\jxo}=\bar y\in\rV$, and \ref{intR1} holds.\vom

{\ubf Case 2:} 
$\jao\in A_1$. 
Then $\muo(A_1)>0$ by the randomness. 
In $\rV$, there is an $\Fs$ set $A'_1\sq A_1$ of the 
same measure, so the Borel set $A_1\bez A'_1$, coded 
in $\rV$, is null, and hence $\jao\in A'_1$. 
Therefore there is, in $\rV$, a perfect set 
$A_2\sq A'_1$, satisfying $\jao\in A_2$ and 
$\muo(A_2)>0$. 
We let $\mu(A)={\muo(A)}/\muo(A_2)$, for any measurable 
$A\sq A_2$, so $\mu$ is a continuous probability 
measure on $P$, and the real $\jao\in P$ is 
$\mu$-random over $\rV$. 
The set $Y_2=\im f{A_2}$ is closed, and by construction 
we have
\ben
\fenu
\itlb{assf1}%
if $x\in Y_2$ then $\mu(\aim fx)=0$ 
(\ie, $f$-preimages of singletons are $\mu$-null). 
\een 

The set $R$ of all rational intervals 
$J\sq\dI$, such that $\mu(\aim f{J\cap Y_2})=0$, 
is at most countable. 
Therefore ${\Ao}=A_2\bez\bigcup_{J\in R}\aim f{J\cap Y_2}$ 
is a closed subset of $A_2$, of the same measure, $f$ 
maps ${\Ao}$ onto the closed set 
${\Yo}=Y_2\bez\bigcup R$, and we have
\ben
\fenu
\atc
\itlb{assf2}%
if $J$ is an open interval in $\dI$ and 
${\Yo}\cap J\ne\pu$ then $\mu(\aim f{{\Yo}\cap J})>0$. 
\een 

\bdf
\lam{hatf}
If $x\in\dI$ then let 
$\haf(x)=\mu(\aim f{{\Yo}\cap[0,x)})$, 
so $\haf:\dI\to\dI$.
\edf

\ble
\lam{xx1}
The map\/ $\haf$ 
is continous, $\ran\haf=\dI$, and 
$\haf$ is strictly increasing, 
except that\/ $\haf(x)=\haf(x')$ in case when\/ 
$x<x'$ belong to\/ $\dI$ and\/ ${\Yo}\cap(x,x')=\pu$.
\ele
\bpf
Let $x<x'$ belong to $\dI$.  
Then $\haf(x)\le\haf(x')$ is clear. 
To prove the strict inequality, note that 
$\haf(x')-\haf(x)=\mu(\aim f{{\Yo}\cap[x,x')})>0$ 
provided ${\Yo}\cap(x,x')\ne\pu$, 
and apply \ref{assf1}, \ref{assf2}. 
\epf

\ble
\lam{xx2}
The superposition map\/ $\baf(a)=\haf(f(a)):{\Ao}\onto\dI$ 
is continuous and measure-preserving in the sense that if\/ 
$X\sq\dI$ is Borel then\/ $\mu(\aim\baf{X})=\la(X)$, 
while if\/ 
$A\sq{{\Ao}}$ is Borel then\/ $\la(\im\baf A)\ge\mu(A)$.
\ele
\bpf
Consider any interval $X=[0,m)$ in $\dI$; $0\le m\le1$. 
By definition, $\haf(x)\in X$ iff 
$\mu(\aim f{{\Yo}\cap[0,x)})<m$.
Therefore the $\haf$-preimage $\aim\haf{X}$ 
is equal 
to $Z=[0,M)$, where $M$ is the largest real in $\dI$ 
satisfying the inequality $\mu(\aim f{{\Yo}\cap[0,M)})\le m$. 
Then clearly $\mu(\aim f{{\Yo}\cap Z})=m$. 

But 
$\aim f{{\Yo}\cap Z}=
\aim f{\aim\haf{X}}=\aim\baf{X}$. 
We conclude that 
$\mu(\aim\baf{X})=\la(X)=m$ for any 
$X=[0,m)$, as above. 
By induction, this implies 
$\mu(\aim\baf{X})=\la(X)$ 
for any Borel $X\sq\dI$, the first claim. 
The second claim follows, since  
$A\sq \aim\baf{\im\baf{A}}$, 
and any analytic set has a Borel subset of 
the same measure. 
\epf

\bcor
[under Case 2]
\lam{xxC} 
The real ${\jyo}=\baf({\jao})=\haf({\jxo})\in\dI$ is 
$\la$-ran\-dom over $\rV$. 
Therefore\/ $\rV[{\jxo}]=\rV[{\jyo}]$ is a Solovay-random 
extension of $\rV$.
\ecor
\bpf 
To prove the second claim, note that $\haf$ is 
\lap{almost} $1-1$ on ${\Yo}$ by Lemma~\ref{xx1}, 
and hence $\rV[{\jxo}]=\rV[{\jyo}]$. 
\epf  

We have another split in cases. 
In $\rV$, let $\cB$ be the family of all Borel sets 
$B\sq \Ao$ such that $\mu(B)>0$ and $F$ is $1-1$ on $B$. 
The set $\cB$ can be empty or not, but anyway there is 
a Borel set $B_0$, equal to a union of $\le\alo$ sets 
in $\cB$, such that $\mu(B'\bez B_0)=0$ for any $B'\in\cB$. 
(If $\cB=\pu$ then $B_0=\pu$ either.) 
We let $\Ai=\Ao\bez B_0$ and $\Yi=\im FB$. 
Thus $\Ai$ is Borel, $\Yi\sq \Yo$ analytic, and 
\ben
\fenu
\atc
\atc
\itlb{assf3}%
if 
$B\sq \Ai$ is Borel and $\mu(B)>0$ then $F$ is {\ubf not} 
1-1 on $B$.
\een

{\ubf Subcase 2a of Case 2:} 
$\jao\in\Ao\bez\Ai$. 
By construction there is a Borel set $B\sq\Ao$ such that 
$\jao\in B$, $\mu(B)>0$, and  $F$ is $1-1$ on $B$. 
Then 
$\jao\in\id11(p,p',\jyo)$ 
for some $p,p'\in\rV$ (codes for $F,B$), 
hence 
${\jao}\in\rV[{\jyo}]=\rV[{\jxo}]$, 
thus (R2) holds.\vom

{\ubf Subcase 2b of Case 2:} 
not Subcase 2a. 
This is the {\ubf key subcase}, 
and it will be considered  
in the two following  sections.

\vyk{
By classical theorems of descriptive set theory, 
the subcase assumption can be reduced 
to the assumption that 
\ben
\fenu
\atc
\atc
\itlb{assf3}%
if $y\in \dI$ then the preimage 
$\aim Fy$ is uncountable. 
\een 
And the goal is to prove \ref{intR2}(b) 
of Thm~\ref{intR} under the assumptions 
\ref{assf1},\ref{assf2},\ref{assf3}. 
}

\parf{The key subcase, measure construction}
\las{XX2}  

Here we prove that $\rV[{\jao}]$ is a random 
extension of $\rV[{\jxo}]$.
First of all, we define, in $\rV[{\jxo}]$, a measure on 
the set $\Om=\aim F{{\jxo}}$, with respect to which 
${\jao}$ itself will be random.
We'll make use of the following lemma which combines 
effects of forcing and the Shoenfield absoluteness 
theorem.

\ble
\lam{xxA}
Let\/ $\vpi(x)$ be a combination of\/ 
$\is11$-formulas and\/ $\ip11$-formulas, by means 
of $\land$, $\lor$, $\neg$, and quantifiers over\/ 
$\om$, and with reals in\/ $\rV$ as parameters. 
If\/ $\vpi(\jyo)$ is true then there is a closed set\/ 
$Y\sq\dI$ of positive measure\/ $\la(Y)>0$, 
coded in\/ $\rV$, containing\/ ${\jyo}$, and 
satisfying\/ $\vpi(y)$ for all\/ $y\in Y$.
\ele
\bpf
The set $\ens{y}{\vpi(y)}$ is measurable, hence, 
it is true in $\rV$ that any Borel $Y_0\sq\dI$ 
of positive measure contains a perfect subset 
$Y\sq Y_0$ still of positive measure $\la(Y)>0$, 
satisfying either (1) $\kaz y\in Y\vpi(y)$ or 
(2) $\kaz y\in Y\,\neg\,\vpi(y)$. 
These formulas are $\ip12$, 
hence absolute by Shoenfield. 
It follows by the randomness of ${\jxo}$ that 
there is a perfect subset $Y\sq\dI$ 
of positive measure, containing ${\jyo}$ and 
satisfying (1) or (2). 
But (2) is impossible because of $\vpi({\jyo})$.
\epf

Now suppose that $B\sq{{\Ai}}$ is a Borel set.

If $X\sq\dI$ then let 
$\og BX=B\cap\aim\baf X=\ens{a\in B}{\baf(a)\in X}$.

In particular, if $x\in\dI$ then 
$\og Bx=\ens{a\in B}{\baf(a)=x}$, a cross-section.

Note that $\mu(B)\le \la(\im\baf B)$
by Lemma~\ref{xx2}, and if $X\sq\im FB$ is Borel 
then $\mu(\og BX)\le \la(X)$. 
If $X\sq\dI$ is Borel then put 
$\la_B(X)=\mu(\og BX)$; 
$\la_B$ is a $\sg$-additive Borel measure 
on $\dI$, concentrated on $\im FB$ and  
satisfying $\la_B(X)\le\la(X)$. 

If $x\in\dI$ then let 
$U_B(x)=\la_B([0,x))=\mu(\og B{[0,x)})$. 
It is important that $U_B:\dI\to\dI$ is   
non-decreasing  
($x<y\imp U_B(x)\le U_B(y)$). 
We'll make use of the following 
collection of classical results 
related to monotone real functions. 

\bpro
[see \eg\ \cite{riss}, Chapters I and II]
\label{ris}
\ben
\renu
\itlb{ris1}%
If\/ $B\sq{{\Ai}}$ is a Borel set then a 
derivative\/  
$U'_B(x)<\iy$ exists for\/ 
$\la$-almost all\/ $x\in\dI\,;$   

\itlb{ris2}%
If\/ $B\sq{{\Ai}}$ is a Borel set and\/ 
$U'_B(x)=0$ for\/ $\la$-almost all\/ $x\in\dI$, 
then\/ $U'_B(x)=0$ for all\/ $x\in\dI\,;$ 

\itlb{ris3}%
if\/ $B_0,B_1,\dots\sq{{\Ai}}$ are 2wise-disjoint 
Borel, and\/ $B=\bigcup_nB_n$, then\/
$U_B(x)=\sum_nU_{B_n}(x)$, $\kaz x$, and\/ 
$U'_B(x)=\sum_nU'_{B_n}(x)$ for\/ $\la$-almost 
all\/ $x\in\dI$.\qed
\een
\epro

\ble
\lam{xx3}
If\/ $C\sq{{\Ai}}$ and\/ $X\sq\im FC$ are Borel sets, 
$\la(C)>0$, and\/ $B=\og CX$, then\/ 
$U_C'(x)=U'_B(x)$ for\/ $\la$-almost all\/ 
$x\in X$. 
\ele
\bpf
Let $A=C\bez B$, so that $X$ and $Y=\im FA$ are 
disjoint sets satisfying $X\cup Y=\im FC$. 
Accordingly we have $U_C(x)=U_B(x)+U_A(x)$ for 
all $x\in \dI$, therefore $U'_C(x)=U'_B(x)+U'_A(x)$ 
for $\la$-almost all $x$ (those in which all three 
derivatives are defined). 
However we have $U'_A(x)=0$ for $\la$-almost all 
$x\in X$; in fact, the equality holds for all 
points $x\in X$ of density 1. 
As required.
\epf

\bdf
\lam{dfom}
Let $\Om=\aim f{{\jxo}}=\aim\baf{{\jyo}}$ 
(a 
closed set, containing ${\jao}$).
\edf

\ble
\lam{xx4}
If\/ $P\sq \Om$ is a Borel set coded 
in\/ $\rV[{\jyo}]$ then there is a Borel set\/ 
$B\sq\dI$, coded in\/ $\rV$ and such that\/ 
$P=\og B {{\jyo}}$.
\ele
\bpf 
There is a Borel set $W\sq\dI\ti{{\Ai}}$, 
coded in $\rV$, such that 
$P=W_{{\jyo}}=\ens{a}{\ang{{\jyo},a}\in W}$ 
(a cross-section). 
Thus $W_{{\jyo}}\sq \Om=\aim\baf{{\jyo}}$. 
By Lemma~\ref{xxA}, there is a Borel set 
$X\sq\dI$ of positive measure $\la(X)>0$, 
coded in $\rV$, containing ${\jyo}$, and such that 
$W_{y}\sq \aim\baf{y}$ holds for all $y\in X$. 
Then 
$B=\ens{a}{\baf(a)\in X\land\ang{\baf(a),a}\in W}$
is a Borel set coded in $\rV$. 
Moreover $W_y=\og B y$ for all $y\in X$ 
by construction, 
in particular, $P=\og B {{\jyo}}$. 
\epF{Lemma}

\bdf
\lam{xxD}
If $P\sq \Om$ is a Borel set coded 
in\/ $\rV[{\jyo}]$ then let $\nu(P)=U_B'({\jyo})$, 
for any $B$ as in the lemma.
It follows from Proposition~\ref{ris}\ref{ris1} 
that $U_B'({\jyo})$ is defined, because ${\jyo}$ is 
random over $\rV$ by the above.
\edf

\ble
\lam{xx6} 
$\nu(P)$ is independent of the choice 
of $B$. 
\ele
\bpf
Suppose that $C\sq {{\Ai}}$ is another Borel set 
satisfying $P=\og C {{\jyo}}$. 
By Lemma~\ref{xxA}, there is a Borel set 
$X\sq\dI$ of positive measure $\la(X)>0$, 
coded in $\rV$, containing ${\jyo}$, and such that 
$\og C {y}=\og B {y}$ holds for all $y\in X$. 
Then $U_B'(y)= U_C'(y)$ for $\la$-almost all $y\in X$ 
by Lemma~\ref{xx3}.
Therefore $U_B'({\jyo})= U_C'({\jyo})$, as ${\jyo}\in X$ 
is random.
\epf

Thus $\nu$ is a well-defined measure on 
Borel sets $P\sq\Om$ in $\rV[{\jyo}]$.

\parf{The key subcase, proof of randomness}
\las{XX3}  

To finalize the proof of Theorem~\ref{intR} 
in Case 2b, we are going to show that 
${\jao}$ is $\nu$-random over\/ $\rV[{\jyo}]$.  
Then it suffices to apply Lemma~\ref{xxl}, 
to transform $\jao$ 
to a \lap{standard} $\la$-random real in $\dI$.
We first of all show that $\nu$ is 
a \lap{good} measure.

\ble
\lam{111} 
In\/ $\rV[{\jyo}]$, 
$\nu$ is a\/ $\sg$-additive continuous  
probability measure on\/ $\Om$. 
\ele
\bpf
(A) 
To prove $\nu(\Om)=1$ take $B={{\Ai}}.$ 
Then $\og B{{\jyo}}=\aim \baf{{\jyo}}=\Om$. 
Lemma~\ref{xx2} implies 
$$
U_B(x)=\la_B([0,x))=\mu(\og B{[0,x)})=
\mu(\aim\baf{[0,x)})=\la([0,x))=x\,,
$$ 
and hence $U'_B(x)=1$ for all $x$.
In particular, $\nu(\Om)=U'_B({\jyo})=1$.\vom

(B) 
Prove $\sg$-additivity.
Lemma~\ref{xx4} reduces this to the 
following claim: 
\rit{if\/ $\sis{C_n}{n<\om}\in\rV$ is a sequence of 
Borel sets\/ $C_n\sq{{\Ai}},$ and\/ 
$(\og{C_k}{{\jyo}})\cap (\og{C_n}{{\jyo}})=\pu$ 
for all\/ $k\ne n$, and\/ $C=\bigcup_nC_n$, 
then\/ $U'_C({\jyo})=\sum_n U'_{C_n}({\jyo})$.} 
By Lemma~\ref{xxA}, there is a Borel set 
$X\sq\dI$ with $\la(X)>0$, 
coded in $\rV$, containing ${\jyo}$, and such that 
$(\og{C_k}{{y}})\cap (\og{C_n}{{y}})=\pu$ 
for all $y\in X$, $k\ne n$. 
The Borel sets $B_n=\og{C_n}X\sq{{\Ai}}$ are 
pairwise disjoint, and the set 
$B=\og CX$ satisfies $B=\bigcup_nB_n$. 

Moreover, we have $U_B(x)=\sum_nU_{B_n}(x)$ 
for all $x$, and\/ 
$U'_B(x)=\sum_nU'_{B_n}(x)$ for $\la$-almost 
all $x\in\dI$ by Proposition~\ref{ris}\ref{ris3}. 
Finally, Lemma~\ref{xx3} implies that 
$U'_B(x)=U'_C(x)$ and $U'_{B_n}(x)=U'_{C_n}(x)$
for all $n$ and $\la$-almost all $x\in X$.
It follows that $U'_C(x)=\sum_nU'_{C_n}(x)$ 
for $\la$-almost 
all $x\in X$, hence, 
$U'_C({\jyo})=\sum_nU'_{C_n}({\jyo})$ 
by the randomness, as required.\vom

(C) 
To prove that $\nu$ is continuous, 
suppose to the contrary that $z_0\in\Om$ and 
$\nu(\ans{z_0})>0$. 
By definition there is a Borel set 
$C\sq{{\Ai}},$ coded in $\rV$ and satisfying 
$\og C{{\jyo}}=\ans{z_0}$ and $U'_C({\jyo})>0$. 
By Lemma~\ref{xxA}, there is a Borel set 
$X\sq\dI$ with $\la(X)>0$, 
coded in $\rV$, containing ${\jyo}$, and such that 
$\og C{{y}}$ is a singleton and $U'_C({y})>0$ 
for all $y\in X$. 
Let $B=\og CX$. 
Then $\og B{{\jyo}}=\ans{z_0}$, 
$\og B{{y}}$ is a singleton for all $y\in X$,  
and $U'_B(y)>0$ 
for $\la$-almost all $y\in X$, by Lemma~\ref{xx3}. 
It follows that $U_B(1)>0$, hence $\mu(B)=U_B(1)>0$.
Moreover, by the singleton condition, the preimage 
$\aim Fy\cap B=\og By$ is a singleton for all 
$y\in\im FB\sq X$. 
But this contradicts the Case 2b assumption. 
\epf

\ble
\lam{xx8} 
${\jao}$ is\/ $\nu$-random over\/ $\rV[{\jyo}]$.
\ele
\bpf
Assume that $P\sq\Om$ is a Borel set, 
coded in $\rV[{\jyo}]$, and $\nu(P)=0$; 
we have to prove that ${\jao}\nin P$. 
By definition there is a Borel set 
$C\sq{{\Ai}},$ coded in $\rV$ and satisfying 
$P=\og C{{\jyo}}$ and $U'_C({\jyo})=0$. 
By Lemma~\ref{xxA}, there is a 
closed (here, this is more suitable than 
Borel) set 
$X\sq\dI$ of positive measure $\la(X)>0$, 
coded in $\rV$, containing ${\jyo}$, and such that 
$U'_C(y)=0$ for all $y\in X$. 

Let $B=\og CX$. 
Then $P=\og B{{\jyo}}$, and $U'_B(y)=0$ 
for $\la$-almost all $y\in X$ by Lemma~\ref{xx3}. 
Note that $\im FB\sq X$, thus $U_B(x)$ is a constant 
inside any open interval disjoint with $X$. 
Thus $U'_B(y)=0$ for all $y\in\dI\bez X$, 
hence overall $U'_B(y)=0$ 
for $\la$-almost all $y\in \dI$. 
This implies $U_B(x)=0$ for all $x\in\dI$ by 
Proposition~\ref{ris}\ref{ris2}.
Therefore $\la_B(\dI)=\mu(B)=0$ by construction. 
We conclude that ${\jao}\nin B$, by the $\mu$-randomness 
of ${\jao}$. 
Then ${\jao}\nin P=\og B{{\jyo}}$, as required. 
\epf

\qeD{Theorem~\ref{intR}}

\bcor
\lam{rl1+}
If\/ $x,y$ are reals in an\/ $\ali$-random extension\/ 
$N=\rL[\sis{a_\xi}{\xi<\omi}]$ of\/ $\rL$, then\/ 
$y$ belongs to a random extension of\/ $\rL[x]$ inside\/ $N$.  
\ecor
\bpf
We have $x\in N_\al=\rL[\sis{a_\xi}{\xi<\al}]$
and $y\in N_\ba$, 
for some $\al<\ba<\omi$. 
The model $N_\al$ is equal to a simple 
extension of $\rL$ by one random real. 
Thus, by Theorem~\ref{intR}, either $N_\al=\rL[x]$ 
or $N_\al$ is a random extension of $\rL[x]$. 
In addition, $N_\ba$ is a random extension 
of $N_\al$. 
This implies the result required. 
\epf

\parf{Proof of the localization lemma}
\las{lok}

\bpf[Lemma~\ref{llok}]
Let $\don$ be the weakest element of any forcing 
considered, and $\dox=\ans\bon\ti x$ be the 
canonical name for any set $x$ in the ground set 
universe $\rV$. 
Let $\raf$ be the random forcing and $\for_\raf$ 
be the associated forcing relation. 

\bcl
\lam{rfc}
If\/ $n\ge2$ and $\vpi(\cdot)$ is a parameter-free\/ 
$\is1n$-formula, resp., $\ip1n$-formula, then the set\/ 
$F_\vpi=\ens{x}{\don\for_\raf\vpi(\dox)}$ is\/ 
$\is1n$, resp., $\ip1n$.
\ecl
\bpf
We make use of a standard Borel coding system for 
subsets of $\dn.$ 
It consists 
of $\ip11$ sets $\kos\sq\dn$ and 
$W_+\yi W_-\sq\bn\ti\bn,$ and 
an assignment $c\mto \bks c\sq\dn$, 
such that 
(1) $\ens{\bks c}{c\in\kos}$ is exactly 
the family of all Borel sets $X\sq\dn,$ and 
(2) if $c\in\kos$ and $x\in\dn$ then 
$x\in\bks c$ iff $W_+(c,x)$ iff $\neg\:W_-(c,x)$. 

To define an associated coding system for Borel maps, 
let $e\mto \sis{(e)_n}{n<\om}$ be a recursive 
homeomorphism $\dn\onto\dnp\om$. 
Let $\kof=\ens{e\in\dn}{\kaz n((e)_n\in\kos}$ --- 
\imar{kof}
codes of Borel maps $f:\dn\to\dn.$ 
If $e\in\kof$ then define a Borel map 
$\bkf e:\dn\to\dn$ so that $\bkf e(x)(n)=1$ iff 
$x\in \bks{(e)_n}$, for all $x\in\dn,$ $n<\om$.

If $\vpi(v_1,\dots,v_k)$ is any formula, 
$e_1,\dots,e_k\in \kof$, and $x\in\bn,$ then let 
$\vpi(e_1,\dots,e_k)[x]$ be the formula 
$\vpi(\bkf{e_1}(x),\dots,\bkf{e_k}(x))$, and let 
$$
\TS
\forc\vpi=\ens{\ang{c,e_1,\dots,e_k}\in\kos\ti\kof^k}
{\mu(\bks c)>0\land 
\bks c\for_\raf \vpi(e_1,\dots,e_k)[\ja]}\,,
$$ 
where $\ja$ is a canonical name for the random real.
We assert the following.
\ben
\fenu
\itlb{rfc*}%
If $\vpi$ is a $\ip11$ formula then $\forc\vpi\in\is12$. 
If $\vpi$ is a $\is1n$ formula, $n\ge2$, 
then $\forc\vpi\in\is1n$. 
If $\vpi$ is a $\ip1n$ formula, $n\ge2$, 
then $\forc\vpi\in\ip1n$. 
\een
This is proved by induction. 
If $\vpi(v)$ is $\ip11$ then $\ang{c,e}\in\forc\vpi$ 
iff the set $X=\ens{x\in B_c}{\neg \:\vpi(\bkf e(x))}$ 
is null, which roughly estimated to be $\is12$ by 
coverings with $\Gd$ sets. 
To pass $\ip1n\to\is1{n+1}$, assume that 
$\vpi(v_1):=\sus v_2\,\psi(v_1,v_2)$, 
$\psi$ is $\ip1n$. 
Then $\ang{c,e_1}\in \forc\vpi$ iff 
$
\sus e_2\in\kof\,(\ang{c,e_1,e_2}\in\forc\psi) 
$.
(We make use of the fact that the random forcing admits 
Borel reading of names.) 
Thus if $\forc\psi$ is $\is1{n+1}$ then so is $\forc\vpi$.
To pass $\is1n\to\ip1{n}$, let 
$\vpi(v)$ be $\is1n$. 
Then 
$$
\ang{c,e}\in \forc{\neg\:\vpi}
\leqv 
\kaz c'\in\kos\,(\bks{c'}\sq \bks{c}\land 
\mu(\bks{c'})>0
\imp\ang{c',e}\nin\forc\psi)\,. 
$$ 
Thus if $\forc\vpi$ is $\is1{n}$ then 
$\forc{\neg\,\vpi}$ is $\ip1{n}$.
This ends the proof of \ref{rfc*}.

Now to prove the claim note that $x\in F_\vpi$ 
iff $\ang{c_0,e_x}\in\forc\vpi$, 
where $c_0\in\kos$ satisfies $\bks{c_0}=\dn,$
while $e_x\in\kof$ is such that $\bkf{e_x}$
is the constant map $\bkf{e_x}(a)=x$, $\kaz a\in\dn$.
\epF{Claim}

To finalize the proof of Lemma~\ref{llok}, 
we define formulas $\vpa(x)$ by induction. 
If $\vpi$ is $\is12$ or $\ip12$ then $\vpa:=\vpi$ 
works by Shoenfield. 
Suppose that $n\ge2$, $\psi(x,y)$ is $\ip1n$, and 
a $\ip1n$-formula $\psa$ is defined, satisfying 
$\psi(x,y)\leqv \rL[x,y]\mo\psa(x,y)$ in 
$N=\rL[\sis{a_\xi}{\xi<\omi}]$ 
(a given $\ali$-random extension). 
We define $\vpa(x)$ to be the formula 
$\don \for_\raf\sus y\,(\rL[\dox,y]\mo\psa(\dox,y))$. 
This is a $\is1{n+1}$-formula by Claim~\ref{rfc}, 
so it remains to show that 
$\vpi(x)\leqv\rL[x]\mo\vpa(x)$ in $N$. 

Assume that $x$ is a real in $N$ satisfying $\vpi(x)$. 
Thus there is a real $y\in N$ satisfying 
$\psi(x,y)$, or equivalently, $\rL[x,y]\mo\psa(x,y)$. 
By Corollary \ref{rl1+}, $y$ belongs to a random 
extension of $\rL[x]$ inside $N$. 
Therefore, as the random forcing is homogeneous, 
it is true in $\rL[x]$ that 
$\don \for_\raf\sus y\,(\rL[\dox,y]\mo\psa(\dox,y))$. 
In other words, $\rL[x]\mo\vpa(x)$.  

To prove the converse, assume that 
$\rL[x]\mo
\big(\don \for_\raf 
\sus y\,(\rL[\dox,y]\mo\psa(\dox,y))\big)$. 
Consider any real $z\in N$ random over $\rL[x]$. 
Then $\sus y\,(\rL[x,y]\mo\psa(x,y))$ holds 
in $\rL[x,z]$, so there is a real $y\in\rL[x,z]$ 
satisfying $\rL[x,y]\mo\psa(x,y)$. 
Then $N\mo\psi(x,y)$ by the choice of $\psa$, 
hence finally $N\mo\vpi(x)$.\vtm  

\epF{Lemma~\ref{llok} and Theorem~\ref{mt'}}

\addcontentsline{toc}{subsection}{\hspace*{4.7ex} References}

\bibliographystyle{plain}
{\small
\bibliography{47,kle}
}

\end{document}